\documentclass[11pt,reqno]{amsart}

\usepackage[italian,english]{babel}

\usepackage{graphicx}
\usepackage{newlfont} 
\usepackage{dsfont}
\usepackage{amssymb}
\usepackage{amsmath}
\usepackage{latexsym}
\usepackage{amsthm}

\usepackage{amsthm}
\usepackage{amsmath}
\usepackage{amssymb}
\usepackage[utf8]{inputenc} 
\usepackage{fancyhdr}
\usepackage{hyperref}
\usepackage{faktor}
\usepackage{mathtools}
\usepackage{thmtools}
\usepackage{tikz}
\usetikzlibrary{matrix,calc}
\usetikzlibrary{cd}
\declaretheorem[name=Theorem,numberwithin=section]{thm}
\declaretheorem[name=Corollary, sibling=thm]{cor}
\declaretheorem[name=Proposition, sibling=thm]{prop}
\declaretheorem[name=Definition, sibling=thm]{defin}
\declaretheorem[name=Lemma, sibling=thm]{lemma}

\theoremstyle{remark}
\newtheorem*{rmk}{Remark}

\newcommand{\N}{\mathbb{N}}
\newcommand{\Q}{\mathbb{Q}}
\newcommand{\C}{\mathbb{C}}
\newcommand{\Z}{\mathbb{Z}}

\newcommand{\op}{\operatorname}

\newcommand{\mc}[1]{\mathcal{[1]}}

\newcommand{\R}{\mathbb{R}}
\newcommand{\PP}{\mathbb{P}}

\newcommand{\img}{\operatorname{Im}}

\let\ForAll\forall
\renewcommand{\forall}{\quad \ForAll}

\newcommand{\eps}{\varepsilon}

\newcommand{\codim}{\operatorname{codim}} 

\newcommand{\ch}{\operatorname{ch}}
\newcommand{\Td}{\operatorname{Td}}
\newcommand{\Stab}{\operatorname{Stab}}
\newcommand{\M}[2]{\bold{M}_{#1}(#2)}
\newcommand{\cM}[2]{\cal{M}_{#1}(#2)}

\newcommand{\Ms}[2]{\bold{M}^s_{#1}(#2)}
\newcommand{\K}[2]{\bold{K}_{#1}(#2)}
\newcommand{\Hcal}{\mathcal{H}}
\newcommand{\dual}[1]{#1\,\check{\vrule height1.3ex width0pt}}
\newcommand{\alb}{\operatorname{alb}}
\newcommand{\Ext}{\operatorname{Ext}}
\newcommand{\Hom}{\operatorname{Hom}}
\newcommand{\ext}{\operatorname{ext}}

\newcommand{\alert}[1]{\par \textcolor{red}{\textbf{\huge{#1}}}\par}
\renewcommand{\alert}[1]{}

\patchcmd{\abstract}{\titlepage}{\thispagestyle{empty}}{}{}
\patchcmd{\endabstract}{\endtitlepage}{}{}{}

\usepackage{appendix}

\usepackage{thm-restate}
\usepackage[capitalise]{cleveref}
\crefname{subsection}{subsection}{Subsection}  
\crefname{prop}{{\bf{Proposition}}}{Proposition}  
\crefname{lemma}{{\bf{Lemma}}}{Lemma}  
\crefname{thm}{{\bf{Theorem}}}{Theorem}  
\crefname{cor}{{\bf{Corollary}}}{Corollary}
\newtheorem*{notation}{Notation}

\author{Giacomo Nanni}
\address{Dipartimento di Matematica \\
	Universit\`a di Bologna\\
	Piazza di Porta San Donato 5\\
	40127 Bologna, Italy}
\email[G.~Nanni]{giacomo.nanni13@unibo.it}

\makeatletter

\title{Classification of Lagrangian planes in generalised Kummer manifolds}

\@namedef{subjclassname@2020}{\textup{2020} Mathematics Subject Classification}
\thanks{The author is part of INdAM
research group “GNSAGA” and was partially supported by it. The author acknowledges the support of the European Union - NextGenerationEU under the National Recovery and Resilience Plan (PNRR) - Mission 4 Education and research - Component 2 From research to business -Investment 1.1 Notice Prin 2022 - DD N. 104 del 2/2/2022, from title "Symplectic varieties: their interplay with Fano manifolds and derived categories", proposal code 2022PEKYBJ – CUP J53D23003840006. During the last part of the reviewing process, the author was also supported by the DFG through the research grant Le 3093/5-1.}

\usepackage[
backend=bibtex,
style=alphabetic,
sorting=nty,
labelalpha=true,
maxbibnames=13,
maxalphanames=13,
url=false,
doi=false,
eprint=false,
isbn=false
]{biblatex}
\AtEveryBibitem{\clearlist{language}}

\usepackage{emptypage}

\newcommand{\hk}{IHS }
\newcommand{\hks}{IHS }

\renewcommand{\Stab}{\operatorname{Stab}^\dagger}
\begin{document}

\cleardoublepage
\maketitle
\begin{abstract}{}
	We prove that the class of a line contained in a Lagrangian plane on a dimension $2n$ hyperk\"ahler manifold $X$ of Kummer type has Beauville-Bogomolov-Fujiki square $-\frac{n+1}{2}$ in $H_2(X,\Z)\cong \dual{(H^2(X,\Z))}$ and  order 2 in the discriminant group of $H^2(X,\Z)$. Vice versa, an extremal primitive ray of the Mori cone verifying these conditions is in fact the class of a line in some Lagrangian plane.
\end{abstract}

\newcommand{\Kumm}{ Kummer }
\newcommand{\KKK}{ $\operatorname{K3}$ }
\newcommand{\KKKs}{ $\operatorname{K3}$}
\newcommand{\Hmuk}{H^*_{alg}}
\newcommand{\vv}{\mathbf{v}}

\section*{Introduction}
The birational geometry of a K3 surface is classically understood in terms of intersection properties of effective classes. A classical result of Kovács \cite{Kovacs} describes the Mori cone $\overline{NE_1}(S)$ (the closed cone of effective curves in $H_2(S,\R)$) of a polarised K3 surface $S$. In particular, he shows that isotropic extremal rays of the Mori cone correspond to elliptic fibrations of the surface, while square (-2) classes are the ones that can be contracted through birational transformations.

In \cite{HT2001}, a program was initiated to achieve a similar description for higher-dimensional irreducible holomorphic symplectic (IHS) manifolds. In \cite[Thesis 1.1]{conjHassTsh}, Hassett and Tschinkel proposed a conjectural description of the Mori cone $\overline{NE_1}(X) \subset H_2(X,\R)$ in terms of the quadratic form induced by the Beauville-Bogomolov-Fujiki (BBF) form $q_X$, via the isomorphism $\dual{H^2(X,\Q)} \cong H_2(X,\Q)$. Specifically, the square of the class of a line in a Lagrangian plane $\mathbb{P}^{\frac{\dim X}{2}} \subset X$ was conjectured to be a deformation-invariant constant $c_X$ of the manifold.

The present work addresses this conjecture. Partial answers were previously given for low-dimensional cases (see \cite{HassTsh, HT, HHT, BJ}). Recently, Bakker \cite{Bakker} resolved the conjecture for IHS manifolds deformation equivalent to the Hilbert scheme of points on a K3 surface, computing $c_X$ explicitly and providing a characterisation for the class of a line contained in a Lagrangian plane.

Our main result is analogous to the one of Bakker for the generalised Kummer deformation type.

\begin{thm} [restated later as \cref{final}]
	Let $(M,h)$ be a polarised holomorphic symplectic manifold of Kummer type and dimension $2n$, with $R \in H_2(M,\Z)$ a primitive generator of an extremal ray of the Mori cone. Then $R$ is the class of a line in a Lagrangian plane if and only if $(R,R) = -\frac{n+1}{2}$ and $2R \in H^2(M,\Z)$.
\end{thm}

Similar to Bakker's work, the result relies on a moduli-theoretic characterisation of extremal Lagrangian planes in certain moduli spaces of Bridgeland-stable objects (\cref{charLag}). This is achieved by studying the birational flops of such planes using techniques introduced by Yoshioka, Bayer, and Macrì (see \cite{YoshPosCon, BayerMac}), which allow an interpretation in terms of wall-crossing of stability conditions.
\subsection*{Acknowledgments}

The author would like to thank his supervisors Giovanni Mongardi and Christian Lehn for their help and suggestions. Particular thanks go to Annalisa Grossi for her invaluable help in reviewing and structuring the paper. The author is also very grateful to Claudio Onorati for his assistance during the revision phase.
\section{Notations and preliminaries}\label{sec:prel}
\newcommand{\Db}{\op{D}^b(S)}
Throughout the paper, $S$ will denote an abelian surface and $\dual{S}$ its dual.
\subsection{Bridgeland stability}
The notion of Bridgeland stability was introduced in \cite{BridgelandOriginal} for a general triangulated category and specialised in \cite{BridgelandK3}to the case of the bounded derived category of coherent sheaves $\Db$ of a surface $S$. We will write $\sigma=(Z_\sigma,\cal P_\sigma)$ for a stability condition of central charge $Z_\sigma$ and slicing $\cal P_\sigma$.

Two objects are said to be S-equivalent if their Harder-Narasimhan (HN) filtrations share the same Jordan-H\"older (JH) factors.
We will denote $\Stab(S)$ the connected component of the stability manifold which includes all stability conditions for which the skyscraper sheaves on $S$ are stable.

With respect to the wall-chamber structure on $\Stab(S)$ for a fixed algebraic Mukai vector $v\in \Hmuk(S,\Z):=H_\Z^0\oplus H_\Z^{1,1}\oplus H_\Z^2$, we will say a stability condition $\sigma$ on a wall $\cal W$ is generic on the wall $\cal W$ if it does not lie on any other wall.

On the Mukai lattice $H^*(S,\Z):=H_\Z^0\oplus H_\Z^{2}\oplus H_\Z^4$ we denote $(v,w):=-v_0w_4+v_2w_2-v_4w_4$ the Mukai product.
By Grothendieck-Riemann-Roch, the Mukai vector map $\vv:=\sqrt{\Td(S)}\ch: K_{num}(S)\rightarrow H^*_{alg}(S,\Z) $ (where $\ch$ denotes the Chern character and $\Td(S)$ the Todd class) provides an anti-isometry of lattices between the numerical Grothendieck group $K_{num}(S)$ equipped with the relative Euler characteristic $\chi$ and the algebraic part of the Mukai lattice $(H^*_{alg}(S,\Z),(,))$. Indeed, for $E,F\in \Db$ we have $(\vv(E),\vv(F))=-\chi(E,F)$.

\subsection{Moduli spaces}
As detailed in \cite{SequivMacr}, fixed a stability condition $\sigma\in \Stab(S)$, a phase $\phi\in\R$ and a Mukai vector $v\in H^*_{alg}(S,\Z)$ one can define moduli stacks $\cM\sigma v$ that parametrise S-equivalence classes of semistable objects of Mukai vector $v$ and phase $\phi$ (which we will always assume to be $\phi=1$).
For generic (meaning, not lying on any wall) stability conditions , these stacks admit coarse projective moduli spaces:
\begin{thm}[\cite{BayerMac} Thm. 1.3,\cite{MMYsomeModuli} Theorem 1.4
	]
	Let $S$ be a smooth projective \KKK or abelian surface and let $v\in  H^*_{alg}(S,\Z)$. If $\sigma\in \Stab (S)$ is generic (does not lie on a wall with respect to $v$) then there is a coarse moduli space $\M{\sigma}{v}$ for $\cM\sigma v$  which is a normal projective irreducible variety, parametrising the S-equivalence classes, with $\mathbb{Q}$-factorial singularities.
\end{thm}
Singularities correspond to strictly semistable sheaves, as simple sheaves are unobstructed and therefore the moduli space of stable objects $\Ms{\sigma}{v}$  (if non-empty) is smooth of the expected dimension $v^2+2$.
If $v$ is primitive, then semistable and stable objects coincide, so that the moduli space $\M{\sigma}{v}$ is smooth of the expected dimension $v^2+2$.

The non-primitive case is described in the following theorem where $S$ is a K3 surface. In the case of $S$ abelian surface an analogous result holds true.
\begin{thm}[\cite{MMP} Theorem 2.15]
	Let $v=mv_0\in H_{alg}^*(S,\Z)$ be a Mukai vector, with  $v_0$ primitive and $m>0$ and let $\sigma\in \Stab(S)$ a generic stability condition with respect to $v$. Then
	\begin{enumerate}
		\item The coarse moduli space $\M{\sigma}{v}$ is non-empty if and only if $v^2+2\geq0$
		\item $\Ms \sigma v\neq \emptyset$ and $\dim \M \sigma v=v^2+2$ if and only if  $m=1$ or $v_0^2>0$ .
	\end{enumerate}
\end{thm}
\begin{rmk}
	In the case of abelian surfaces, as there are no rigid objects (after \cite[Lemma 15.1]{BridgelandK3}), so the inequality $v^2+2\geq 0$ is replaced by $v^2\geq 0$.
\end{rmk}

\begin{prop}[from \cite{BayerMac} Lemma 7.2]\label{primIsotrStab}
	Let $v\in H^*(S,\Z)$ be of the form $v=mv_0$, with $v_0^2=0$ and $m>1$. Then $$\M{\sigma}{v}\cong\operatorname{Sym}^m(\M \sigma {v_0})$$
\end{prop}

\subsection{ Kummer \hks inside moduli spaces}\label{subsec:KummerConstruction}
Given an abelian surface $S$, a Mukai vector $v$ with $v^2\geq 2$ and a generic stability condition $\sigma\in \Stab(S)$, then the Albanese map of $\M \sigma v$ can be described as follows (as explained in \cite{YoshPosCon}).

Fix an $E_0\in \cM \sigma v (\C)$. Then define
\begin{alignat*}3
	\alb: & \M\sigma v & \longrightarrow & S\times \dual S                                \\
	      & E          & \mapsto         & (\det\Phi_{\boldsymbol{P}}(E-E_0),\det(E-E_0))
\end{alignat*}
Here, $\boldsymbol{P}$ is the Poincaré bundle on $S\times \dual{S}$ defined by the $\dual{S}$ parametrising line bundles 
on $S$, and $\Phi_{\boldsymbol{P}}$ denotes the corresponding Fourier-Mukai transform.

\begin{defin}
	Let  $v\in H_{alg}^*(S,\Z)$ be a primitive Mukai vector with  $v^2\geq 6$ and $\sigma\in\Stab(S)$ generic. We denote $\K \sigma v$ the fibre  $\alb^{-1}(0)$.
\end{defin}
\begin{thm}[\cite{YoshPosCon} Theorem 1.13]
	Let $v,\sigma$ be as in the previous definition. Then $\K \sigma v$ is an irreducible symplectic manifold of $\dim \K \sigma v=v^2-2$ which is deformation equivalent to a generalised Kummer manifold.
\end{thm}
\subsection{Walls, lattices and birational models}

Fix the Mukai vector $v$ and take a wall $\cal W\subset \Stab(S)$. Then the {associated lattice} is defined as $$\cal H_{\cal W}:=\{w\in \Hmuk(S,\Z), Z(w)/Z(v)\in\R \forall (Z,\cal P)\in\cal W \}.$$ Conversely, for any signature (1,1) sublattice $\cal H\subset \Hmuk(S,\Z)$ containing $v$ we can define the potential wall $$\cal W_{\cal H}:=\{(Z,\cal P)\in \Stab(S), Z(w)/Z(v)\in \R \forall w\in \cal H\}.$$

If $\sigma_+,\sigma_-\in \Stab(S)$ are stability conditions in two adjacent chambers separated by a wall $\cal W\subset\Stab(S)$, then the two moduli spaces $\M {\sigma_+} v,\M {\sigma_-} v$ are birational to each other (see \cite{MMP,YoshPosCon}).

More explicitly, one can construct a birational transformation $\M {\sigma_-} v\dashrightarrow\M{\sigma_+} v$ that identifies objects that are both $\sigma_+$ and $\sigma_-$ stable. This birational map is said to be obtained by wall-crossing.

As a consequence of the construction (see \cite[Prop. 5.1]{MMP}), if a curve in $\M{\sigma_-} v$ is contracted by the wall crossing map, the Mukai vectors of the $\sigma_+$ JH factors of a generic point of the curve are contained in the lattice $\cal H_{\cal W}$ associated to the wall $\cal W$.
This allows to classify wall-crossing transformations in terms of numerical properties of $\cal H_{\cal W}$ (see \cite[Prop. 5.15]{MMP}).

Moreover, we have a wall-crossing interpretation of all IHS birational models of moduli spaces:
\begin{thm}[{\cite[Corollary 3.33]{YoshPosCon}}] \label{mainMMP}
	Let $S$ be an abelian surface. Fix $\sigma\in \Stab(S)$ generic with respect to a primitive positive Mukai vector $v\in H^*_{alg}(S,\Z)$.
	Every K-trivial smooth birational model of $\K\sigma v$ (resp. $\M \sigma v$) appears as a moduli space $\K{\hat\sigma} v$ (resp. $\M {\hat\sigma} v$) of Bridgeland stable objects for some generic $ \hat\sigma\in \Stab(S)$.
\end{thm}
\subsection{Mukai homomorphism}\label{MukaiHom}
Given a (quasi)-universal object $E\in \op{D}^b(S\times \M {\sigma} v)$, we will denote $\theta$ the Mukai homomorphism given by the composition
\begin{equation*}
	\theta:	v^\perp\xrightarrow{-\vv^{-1}} K_{num}(S)\xrightarrow{\Phi_E} K_{num}(\M {\sigma} v)\xrightarrow{\det} H^2_{alg}(\M \sigma v,\Z)\rightarrow H^2_{alg}(\K \sigma v,\Z)
\end{equation*}

where $\Phi_E:K_{num}(S)\rightarrow K_{num}(\M {\sigma} v)$ is the K-theoretic Fourier-Mukai transform.

The dual of its inverse will be denoted $\dual\theta:H^*_{alg}(S,\Z)\rightarrow H^{alg}_2(\K \sigma v)$.

According to the work in \cite[Theorem 4.9]{Wieneck} 
it is possible to extend this setting to all \hk manifolds of generalised Kummer type. For such a manifold $M$, he defines a lattice together with a Hodge structure $\tilde{\Lambda}(M)$ (called the Mukai lattice of $M$) and a monodromy invariant embedding
\begin{equation*}
	H^2(M,\Z)\hookrightarrow\tilde{\Lambda}(M)
\end{equation*}
which for moduli spaces coincides with $\theta^{-1}$. We will still denote as $\dual\theta$ the dual of this embedding.

\subsection{Lagrangian planes}\label{subsec:Lagrangian}
Let $X$ be an \hk manifold of dimension $2n$.
\begin{defin}
	A {Lagrangian plane} in $X$ is a maximal isotropic submanifold isomorphic to a projective space.
	Equivalently, it is a submanifold isomorphic to $\PP^n$.
\end{defin}

Lagrangian planes are rigid in the following sense:
\begin{prop}[{\cite[Lemma 11]{Bakker}}]\label{rigidity}
	Let $\PP\subset X$ be a Lagrangian plane in a holomorphic symplectic variety. Then no curve $C\subset \PP$ deforms outside of $\PP$.
\end{prop}
We say that a Lagrangian plane is {extremal} if the class of its lines lies on the boundary of the Mori cone.
\begin{rmk}
	On a projective IHS manifold $X$, an extremal class $R\in \overline{NE_1}(X)$ with negative BBF square is contractible in the sense of MMP as a consequence of the cone theorem (stated as in \cite[Theorem 3.7]{KollarMori}. To apply the theorem, one needs to choose an effective $\R$-divisor $D$ for which the intersection product $R.D$ is negative. Denote $l\in H^2(X,\Z)$ the dual class of $R$, so that $(l,E)=R.E$ for each $E\in H^2(X,\Z)$. The extremality of $R$ implies that $R^\perp\subset H^2(X,\R)$ is a face of the cone $\op{Nef(X)}=\overline{\op{Amp(X)}}$ closure of the ample cone. Let $H$ be a polarisation of $X$. Then there is $\eps\in\R,\eps>0$ small enough such that $(l+\eps H).R=l.R+\eps H.R=(R,R)+\eps H.R<0$. Therefore $D:=l+\eps H$ is the needed divisor.
\end{rmk}

\newcommand{\pit}{\text{primitive isotropic }}

\subsection{$\PP$-type lattices}\label{sec:p-type}
Let $S$ be an abelian surface and let $A\in \Db$ be an object. Since an abelian surface doesn't admit rigid objects (see \cite[Lemma 15.1]{BridgelandK3}), Serre duality implies that isotropic Mukai vectors are those of minimal square. Moreover, if $A$ is simple (meaning, it doesn't admit non-trivial
subobjects) then its Mukai vector must satisfy $\vv(A)^2=0$ and it must be primitive (as a consequence of \cref{primIsotrStab}).

Notice that this numerical condition is sufficient to grant the stability of a semistable object:
\begin{prop}\label{JHminSq}
	Let $B\in D^b(X)$ be a $\sigma$-semistable object (for some stability condition $\sigma$) having a \pit Mukai vector. Then  $B$ is simple.
	\proof
	Suppose $B$ is not simple.
	and has primitive isotropic Mukai vector $b$. Then, letting $a_1,...,a_m$ be the Mukai vectors of its J-H factors, there would be a map $$\prod_i \Ms{\sigma}{a_i}\rightarrow \M{\sigma}{b}, A_1,...,A_m\mapsto \bigoplus A_i$$ which locally in the image has a section (taking the J-H factors, by the same technique of \cite[Lemma 6.5]{BayerMac})
	. The space on the left must have a smaller dimension than the one on the right, which means $ \sum_i a_i^2 +2m\leq b^2+2$. Since $m\geq2$, this implies $b^2\geq a_i^2$ for each $i$. So by the minimality we get $a_i^2=0$.
	Applying \cref{primIsotrStab}, we conclude that each stable factor must be \pit.
	Moreover, the same inequality shows $m=1$, proving that $B$ is itself its unique stable factor.\endproof
\end{prop}


Once \cref{JHminSq} is established, we can recover the following analogue of $\text{\cite[Proposition 6.3]{MMP}}$ for abelian surfaces
\begin{cor}\label{latticeMMP}
	Let $\mathcal{W}$ be the potential wall associated with a \alert{?}primitive hyperbolic lattice $\Hcal$, and $\sigma\in \mathcal{W}$ generic.
	Then one of the following holds:
	\begin{itemize}
		\item There are no \pit classes in $\Hcal$
		\item There are exactly  two (linearly independent) \pit classes $a,b\in\Hcal$ and there exist $A,B \in\mathcal{P}_\sigma(1)$ with $\vv(A)=a,\vv(B)=b$.
	\end{itemize}
	\proof
	If $a\in\Hcal$ is a \pit class, $(a,a)\geq 0$ so the moduli space $\M{\sigma}{a}$ is non-empty. In other words, there exists a $\sigma$-semistable object $A\in\mathcal{P}_\sigma(1)$ with $\vv(A)=a$ that, by \cref{JHminSq}, is simple.

	The boundary of the positive cone of the realification $\mathcal{H}_\R$ is composed of two lines, which either are both rational or both irrational. In the first case we have two (linearly independent) primitive isotropic vectors $a,b\in \mathcal H$, in the second none.

	\endproof

\end{cor}

We give the following definitions:
\begin{defin}
	Let $S$ be an abelian surface and fix a Mukai vector $v\in H^*_{alg}(S,\Z)$.
	A rank 2 saturated sublattice $\cal H\subset H^*_{alg}(S,\Z)$ of the Mukai lattice is said to be:
	\begin{itemize}
		\item		pointed if $v\in \cal H$
		\item of $\PP$-type if
		      $$\frac{v^2}{2}=\min |(a,v)|$$
		      where the minimum is taken over primitive isotropic classes $a\in \cal H$.
	\end{itemize}
\end{defin}
Using the previous result we can prove that:
\begin{prop}\label{decomposition}\label{bastaLaScomposizione}
	A  rank 2 pointed sublattice $\cal H $ is of $\PP$-type if and only if for a generic $\sigma\in \cal W_{\cal H}$ there exist two simple objects $A,B\in \cal P_\sigma(1)$ with \pit Mukai vectors such that

	\begin{itemize}
		\item $v= \vv(A)+\vv(B)$
		\item $(\vv(A),v)=\frac{v^2}{2}$
	\end{itemize}
	\proof
	$(\implies)$
	Let $a\in\Hcal$ be one of the \pit classes realising the minimum $(a,v)=\frac{v^2}{2} $ in the definition of $\PP$ type.
	Then $b:=v-a$ is \pit:
	\begin{equation*}
		b^2=v^2+a^2-2(v,a)=a^2=0
	\end{equation*}
	and $b$ also realise the minimum
	\begin{equation*}
		(b,v)=v^2-(a,v)=\frac{v^2}{2}.
	\end{equation*}
	This also implies that $b$ is primitive: otherwise $b=m\hat{b}$ would give $\frac{v^2}{2}=m(\hat{b},v)$ and $(\hat{b},v)<\frac{v^2}{2}$ which contradicts the fact that $\frac{v^2}{2}$ is a minimum.
	By \cref{latticeMMP} we get that there exist two simple objects $A,B\in \mathcal{P}_{\sigma_0}(1)$ having linearly independent  \pit Mukai vectors $a,b$.\\

	$(\impliedby)$
	Clearly $a=\vv(A),b=\vv(B)$ are \pit and by \cref{latticeMMP} they are the only two coming as \pit Mukai vectors of simple objects of phase 1.	So for any other \pit class $\alpha$ we can write $\alpha=xa+yb$ for $x,y\in \N$ not both zero. Therefore $|(\alpha,v)|=|(x+y)|\frac{v^2}{2}\geq\frac{v^2}{2}$
	\endproof
\end{prop}
\section{Results}
\subsection{Construction of Lagrangian planes}\label{sec:construct}
Given a $\PP$-type sublattice, we will construct a $\PP^n$-bundle parametrising the H-N filtrations compatible with the decomposition of $v$ coming from \cref{bastaLaScomposizione}. Some of its fibres will be Lagrangian planes in the Kummer manifold $\K{\sigma}{v}$.\\

We define a {partition} of $v\in H^*_{alg}(S,\Z)$ as a set of vectors $\{a_1,...,a_n\}$ summing to $v$, and we denote it $P=[v=\sum_i a_i]$.
For a partition $P=[v=\sum_i a_i]$ and a stability condition $\gamma\in \Stab(X)$ we denote $M_{\gamma,P}\subset \M{\sigma}{v}$ the closed substratum of S-equivalence classes having J-H partition $P$ with respect to $\gamma$.


Assume for simplicity that there exists a universal family $\mathcal{E}$ on $M_{\gamma,P}$. Then by \cite[Theorem 4.3]{MMP} there exists a system of maps \begin{center}
	\begin{tikzcd}
		0=\mathcal E^0\arrow[r]&\mathcal E^1\arrow[r]&...\arrow[r]&\mathcal E^n=\mathcal{E}
	\end{tikzcd}
\end{center} such that on an open set $U\subset M_{\gamma,P}$ its restriction to $u\in U$ is the $\gamma$-HN filtration of $\mathcal E_{u}$. The quotients $\mathcal A_i:=\faktor{\mathcal{E}_{i+1}}{\mathcal{E}_i}$ are families of $\gamma$-semistable objects of Mukai vectors $a_i$. Therefore, there is (similarly to the proof of \cite[Lemma 6.5]{MMP}) a rational map
\newcommand{\JH}{\mathcal{JH}}
\begin{equation*}
	\JH:M_{\gamma,P}\dashrightarrow \M{\gamma}{a_1}\times...\times\M{\gamma}{a_n}
\end{equation*}
which associates to a point in the open set of definition the tuple of its HN semistable factors. 

When the partition has length 2, the fibre over a pair $(A,B)$ of $\gamma$-stable complexes can be identified with $\PP:=\PP(\Ext^1(A,B))$ as follows. Denote $p_S:\PP\times S\rightarrow S$ the second projection, then the HN filtration of any complex in the fibre over $(A,B)$  arises taking the fibre of the distinguished triangle
\begin{equation*}
	(p_S^*B)(1)\rightarrow E\rightarrow p_S^*A\rightarrow(p_S^*B)(1)[1]
\end{equation*}
on points $x\in \PP$
(see \cite[Lemma 15]{Bakker}).\\

We produce some extremal Lagrangian planes on $\K\sigma v$ as intersection of the bundle generated by the map $\JH$ for the partition associated to a $\PP$-type lattice:
\begin{prop}\label{construction}
	Let $\Hcal\subset H^*_{alg}(S,\Z)$ be a $\PP$-type sublattice, then there are Lagrangian planes in $\K \sigma v$.
	In particular, for generic  $\sigma_-$ in the adjecent chamber on the other side of $\mathcal{W}_\Hcal$ and $P=[v=a+b]$ the corresponding partition, there is an open set $U\subset M_P$ such that $\K \sigma v\cap U$ is a union of isolated Lagrangian planes.
	\proof
	By \cref{decomposition}, we have the decomposition $P$ and two simple objects $A,B\in\mathcal{P}_{\sigma_-}(1)$ with Mukai vectors $a,b$, respectively.

	By the previous discussion, we get an open dense set $U\subset M_{\sigma_-,P}$  with a map ${\mathcal{JH}:U\rightarrow \Ms{\sigma_-}{a}\times\Ms{\sigma_-}{b}}$, which is a fibre bundle with fibre $\PP^n$, where $n+1:=(a,b)$.

	It is easy to verify that $\dim \K \sigma v=2n$:
	\begin{equation*}
		2n=2(a,b)-2=2(a,v)-2(a,a)-2=v^2-2=\dim \K \sigma v
	\end{equation*}
	We need to prove that:
	\begin{enumerate}
		\item Fibres are contained in a translate of $\K{\sigma}{v}$.
		\item The intersection of $\K{\sigma}{v}$ with the base is transverse. By intersection with the base, we mean the intersection with the image of the base through the section $\Ms{\sigma_-}{a}\times\Ms{\sigma_-}{b}\rightarrow \M{\sigma}{v}$ given by taking the direct sum.
	\end{enumerate}

	For the first claim, take a fibre $\PP$, which we know to be isomorphic to $\PP^n$.
	Then it would induce a rational curve in $S\times \dual S$ via the Albanese map as $
		\PP\subset\M{\sigma}{v}\rightarrow S\times\dual{S}
	$.
	But since there are no rational curves on an abelian variety, $\PP$ is contained in a fibre of the Albanese map, which by definition is a translate of $\K{\sigma}{v}$.

	The second statement follows from the first. Since each fibre is contained in one of the fibres of the Albanese map, it is either contained in $\K\sigma v$ or disjoint.
	If the planes were not isolated then we would get a deformation of the Lagrangian plane inside the Kummer. But this is not possible because of \cref{rigidity}
	\endproof
\end{prop}
\begin{rmk}
	Since all the lagrangian planes in $\K \sigma v\cap U$ are fibers of the same bundle, the class of a line contained in one of them is the same as the class of a line contained in any other.
\end{rmk}
On planes arising from this construction, we can calculate the class of a line by means of the projection formula. The argument is the same as in \cite[discussion after Lemma 16]{Bakker} but we reproduce it for clarity. In the following lemma, $R\in H_2(S,\Z)$ is the class of a line on the Lagrangian plane $\PP$ obtained by the $\PP$-type lattice $\Hcal$ and $a,b\in \Hcal$ are the two \pit classes realising the minimum, as in the construction of \cref{construction}.
\begin{prop}
	With the notation above, $R=\pm \dual\theta(a)$ where $\theta$ is the Mukai homomorphism.
	Moreover, 
	$\PP$ is extremal.
	\proof
	Denote $(A,B)$ the pair of  J-H factors of objects in $\PP$ for a generic condition $\sigma_0\in \cal W_{\cal H}$.
	Without loss of generality, we can assume $v(A)=a$,$v(B)=b$ and $\phi(A)< \phi(B)$.
	Take $\mathcal{C}$ a curve in $\M{\sigma}{v}$. Then \begin{equation*}
		\theta(w).\mathcal{C}=(w,\vv(\Phi_E(\mathcal{O}_C)))
	\end{equation*}
	To prove the claim it is enough to assume $w\in v^\perp$. If $R$ is the class of a certain $\PP^1\subset\M{\sigma}{v}$ then
	\begin{equation*}
		\theta(w).R=(w,\vv({p_S}_*E_{|\PP^1\times S}))
	\end{equation*}
	By projection formula, we can write (in K-theory):
	\begin{equation*}
		{p_S}_*E_{|\PP^1\times S}=	{p_S}_*((p_S^*A)(1)_{|\PP^1\times S})+	{p_S}_*p_S^*B_{|\PP^1\times S}=A\cdot {p_S}_*p_{\PP^1}^*\mathcal{O}_{\PP^1}(1)+B
	\end{equation*}
	Since by standard computation ${p_S}_*p_{\PP^1}^*\mathcal{O}_{\PP^1}(1)=2$, we conclude ${p_S}_*E_{|\PP^1\times S}=2A+B$.
	Substituting we get:
	\begin{equation*}
		\theta(w).R=(w,\vv(2A+B))=(w,2a+b)=(w,v+a)=(w,a)
	\end{equation*}
	On the other side of the wall, we have $\phi(A)> \phi(B)$, so we get $R=\dual\theta(b)=\dual{\theta}(-a)=-\dual{\theta}(a)$. As crossing $\mathcal{W}_\Hcal$, $R$ changes sign $\PP$ is extremal. \endproof
\end{prop}

\subsection{Characterisation of extremal Lagrangian planes}\label{sec:char}
In this section, we show that every extremal Lagrangian plane can be obtained by the construction from the previous section. We consider the contraction of the Lagrangian plane as wall-crossing (in the sense of \cite{YoshPosCon}) to show that the associated lattice is of $\PP$-type and apply \cref{construction}.




Before entering the proof of the main result of the section (\cref{charLag}), we recall that  there exists a notion of JH filtration in families for families of semistable objects parametrised by a curve.

\begin{prop}[{\cite[Lemma 3.9]{BayerMac}} ]\label{prop:JHfiltrationFamilies}
	Let $C$ be an integral curve. Let $\cal C\in D(S\times C)$ be a family of $\sigma_0$-semistable objects, generically S-equivalent. There exists a filtration :
	$$0=\cal C_0\rightarrow\cal C_1\rightarrow...\rightarrow \cal C_k=\cal C$$
	such that restricting to a generic point $c\in C$ gives the JH filtration of $\cal C_{|\{c\}\times S}$. Moreover, there exist line bundles $\cal L_i$ on $C$ and $\sigma_0$-stable objects $F_i\in D(S)$ such that  $\faktor{\cal C_i}{\cal C_{i-1}}\cong A_i\boxtimes \cal L_i$.
\end{prop}

We will also need this technical fact on deformations of families of semistable objects, the proof of which is delayed to the end of the section.


\begin{lemma}\label{prop:defFamilies}
	Let $\cal T_i \in D(S\times T_i)$ be deformations of $A_i$ over a base $T_i$, such that $\cal T_i\cong \cal T_j$ if $A_i\cong A_j$. Denote $J$ the diagram with objects $T_i$ and morphisms the ones induced by the isomorphisms between $\cal T_i\cong\cal T_j$'s for $i<j$. Then there exists a deformation $\cal F\in D(S\times C\times T)$ of $\cal C$ over the limit $T=\varprojlim_J T_\bullet$ (possibly restricting to open subsets of $T_i$)
	such that the JH-factors of $\cal F_{|S\times \{c\}\times \{t\}}$ are the ${\cal T_i}_{|S\times \{t\}}$.
\end{lemma}
\begin{rmk}
	We recall that $T=\varprojlim_JT_\bullet$ is defined as the closed subvariety of $\bigtimes_i T_i$ given by the equations $\phi_{i,j}(t_i)=t_j$ where $\phi_{i,j}$ are the isomorphisms between $ T_i\cong T_j$.
\end{rmk}

\begin{prop}\label{charLag}
	Let $v\in\tilde{H}_{alg}(S,\Z)$ be primitive and $v^2>0$, and let $\sigma$ be a generic stability condition with respect to $v$.
	If $\PP\subset \K \sigma v$ is an extremal Lagrangian plane, then there exists a $\PP$-type sublattice $\Hcal\subset\tilde{H}_{alg}(S,\Z)$. In particular, there is some stability condition $\sigma_-$ such that $\PP$ is one of the component of $M_ {\sigma_- ,P}\cap \K \sigma v$ where $P$ is the partition associated to $\Hcal$.
	\proof
	The contraction $\phi:\M{\sigma}{v}\rightarrow \bar M$ of $\PP$ is obtained as wall-crossing by \cref{mainMMP}, meaning that a curve is contracted if and only if objects parametrised by such curves are (generically) S-equivalent for some stability condition $\sigma_0$ on an adjacent wall (see \cite[Theorem 2.19]{MMP}). Therefore the Jordan-Holder stable factors of objects in $\PP$ are all of the same form $A_1,...,A_k$. In particular, considering a stability condition $\sigma_-$ on the other side of the wall, $\PP$ is contained in the stratum $M_{\sigma_-,P}$ associated to the partition of the Mukai vector $P=[v=\sum_i a_i]$ where $a_i=v(A_i)$ are the Mukai vectors of the $A_i$.

	Take an open subset of a curve $C\subset \PP$  . From \cite[Lemma 3.9]{BayerMac}, up to restricting to some open subset of $C$, the restriction  $\cal C$ of the universal family of $\M \sigma v$ to $C$ admits a $\sigma_0$-JH filtration given by $$0\rightarrow \cal C_1\rightarrow...\rightarrow \cal C_k$$ where the families of stable factors $\cal A_i $ are constant families of factor $A_i$. Let $\cal T_i$ be the restriction of the universal families of $\M{\sigma_0}{a_i}$ to a neighborhood of $A_i\in \M{\sigma_0}{a_i}$. By \cref{prop:defFamilies} we get a deformation $\cal F\in D(S\times C\times T)$ of $\cal C$ over an open set of $T=\varprojlim_J T_\bullet$ and therefore a map $T\times C\rightarrow \M{\sigma}{v}$. Moreover by construction, each curve image of ${t}\times C, t\in T$ parametrises objects generically of the $\sigma_0$ S-equivalence class parametrised by $t$. Therefore, we get the following commutative diagram
	\begin{center}
		\begin{tikzcd}
			T\times C \arrow[r] \arrow[d] & M_P\arrow[d, "\phi"]\\
			T \arrow[r] & \bar M
		\end{tikzcd}
	\end{center} where the map $T\rightarrow \bar M$ is injective.
	In particular the image of $T\times C$ in $M_P$ gives a deformation of the curve $C$ inside the stratum $M_P$ over the base $T$ such that no other fiber than $C$ is contained in $\PP$. From \cite[Lemma 11]{Bakker} the curve $C$ cannot deform in $\K{\sigma}{v}$ outside $\PP$. This implies that $\dim T\leq \codim \K{\sigma}{v}=4$. The dimension of $T$ is the sum of the terms $a_i^2+2$ without repetitions for isomorphic $A_i$. One can write $I=\{\min\{i, A_i\cong A_j\}, j=1,...,k\}$ and $\dim T=\sum_{i\in I} (a_i^2+2)$, and since $a_i^2\geq 0$ we have $4\geq \dim T\geq 2\#I$, where $\#I$ denotes the cardinality of $I$.  Since $\PP$ is not a point, there are at least two factor $A_i$ and at least two of them are non-isomorphic (as otherwise $v$ would not have been primitive). Therefore $\#I=2$. Denoting $I=\{b_1,b_2\}$ we can then rewrite $4+b_1^2+b^2_2\leq 4$ so $b_1^2=b^2_2=0$. Therefore, we have that shown that each $a_i$ is isotropic and either $a_i=b_1$ or $a_i=b_2$. In other words, the partition $P$ is of the form $v=[b_1+..+b_1+b_2+...+b_2]$, a priori with possible repetition. However, \cite[Lemma 4.3.4 (1)]{YoshFourier} implies that no repetition can occur: if $a_i=a_j,a_j^2=0$ then also $(a_i,a_j)=0$, and neither of the two possible cases applies. So there must be exactly two isomorphic Mukai vectors $a_i$, meaning that the partition is of the form $P=[v=a_1+a_2]$.

	The $\PP$-type lattice $\cal H$ is the one generated by $a_1,a_2$ and by construction $\PP\subset M_{\sigma_-,P}\cap \K{\sigma}{v}$.

	\endproof
\end{prop}

It remains to prove \cref{prop:defFamilies}. Before entering the proof, we fix some notation.

\begin{notation}
	We will denote $T_{\leq m}=\varprojlim_{J_{\leq m}} T_\bullet$ where $J_{\leq m}$ is the restriction of $J$ to the first $m$ indices.
	We will abusively denote by $\cal T_i$ also the pullback of $\cal T_i$ to $S\times T$ or $S\times T_{\leq m}$ along the projections $T\rightarrow T_i$ or $T_{\leq m}\rightarrow T_i$. We will do similarly for any other family defined on $T_i$. We will always allow to restrict to open subsets of the $T_i$.

	To lighten the notation we will write $\cal T_r^t:={\cal T_r}_{|S\times \{t\}}$ and $\cal F_m^{c,t}:={\cal F_m}_{|S\times\{c\}\times \{t\}}$.
	Each time there is no risk of confusion, we will denote the projection for a product $X\times Y$ as $\pi_X,\pi_Y$. Given an object $G\in D(S\times X)$ we will denote $\cal G^Y\in D(S\times X\times Y)$ the constant family over $Y$ given by $\cal G^Y=\pi^*_{S\times X}\cal G$.
\end{notation}
\begin{proof}[Proof of \cref{prop:defFamilies}]
	The proof proceeds by induction: we show that for every $m\leq k$ there exists a family $\cal F_m\in D(S\times C\times T_{\leq m})$ over $T_{\leq m}$ such that
	\begin{itemize}
		\item $\cal F_m$ is a deformation of $\cal C_m$ over $T_{\leq m}$.
		\item the JH factors of $\cal F_m^{c,t}$ are the ${\cal T_i}^t$, $i\leq m$.
		\item the dimension of the group $\Ext^1(\cal T_{m+1}^t,\cal F_m^{c,t})$ does not depend on $t$ and $c$.
	\end{itemize}
	For $m=1$, take for $\cal F_1:=\cal T_1^C$. This verifies that $${\cal F_1}_{|S\times C\times \{0\}}= {(\pi^*_{S\times T_1}\cal T_1)}_{|S\times C\times \{0\}}=\pi^*_{S\times \{0\}}{\cal T_1}_{|S\times 0}= \pi^*_{S\times \{0\}}A_1=\cal A_1=\cal C_1$$
	and the last condition on the $\Ext^1$ groups is verified by virtue of the hypothesis $A_i\cong A_j \iff \cal{T}_i\cong \cal{T}_j$: indeed ${\cal F_1}^{c,t}={\cal T_1}^t$ and therefore $\ext^0({\cal T_2}^t,{\cal T_1}^t)$ is 1 if $A_1\cong A_r$, 0 otherwise, $\ext^2=\ext^0$ is therefore constant by Serre duality and $\ext^1=(a_1,a_r)-2$ or respectively $(a_1,a_2)$, in both cases it does not depend on $t$.

	We therefore assume to have $\cal F_{m-1}$ as above. We start from the triangle $$\cal C_{m-1}\rightarrow \cal C_{m}\rightarrow \cal A_ m.$$ This corresponds to a morphism $\delta^0\in \cal Hom (\cal{A}_m,\cal C_{m-1}[1])=\cal Ext^1 (\cal{A}_m,\cal C_{m-1})$.
	Since by inductive hypothesis the dimension of the fibers does not change as $t\in T_{\leq m}$ varies, $\cal Ext^1 (\cal T_m^C,\cal F_{m-1})$ is locally free on $T_{\leq m}$ and therefore we can (up to restricting $T$) extend $\delta^0$ to a morphism $\delta\in \cal Ext^1 (\cal T_m^C,\cal F_{m-1})$ trivially. We define $\cal F_m:=\op{Cone}\delta$. By construction $${\cal F_m}_{|S\times C\times \{0\}}=\op{Cone}(\delta_{|S\times C\times \{0\}})=\op{Cone}(\delta^0)=\cal C_m$$ and the JH factors for $\cal F_m^{c,t}$ are by construction $\cal T_i^t, i\leq m$. It remains to prove the last condition on the $\ext^1$.

	We consider the triangle $\cal F_{m-1}\rightarrow \cal F_m\rightarrow \cal T_m^C$. Restricting we have triangles $ \cal F_{m-1}^{c,t}\rightarrow \cal F_m^{c,t}\rightarrow \cal T_m^t$ Applying the functor $Hom(\cal{T}_r^t, -)$ we obtain:\begin{align*}
		\dots \rightarrow Hom({\cal T}_r^t, {\cal T}_m^t) \xrightarrow{d_{c,t}} Hom({\cal T}_r^t, {\cal F}_{m-1}^{c,t}[1]) \xrightarrow{f_{t,c}} Hom({\cal T}_r^t, {\cal F}_{m}^{c,t}[1]) \\
		\xrightarrow{g_{t,c}} Hom({\cal T}_r^t, {\cal T}_m^t[1]) \xrightarrow{h_{c,t}} Hom({\cal T}_r^t, {\cal F}_{m-1}^{c,t}[2]) \rightarrow \dots
	\end{align*}

	\renewcommand{\img}{\op{img}}
	From this, since the sequence is exact, we obtain that
	\begin{align*}
		\ext^1 ({\cal T}_r^t, {\cal F}_{m}^{c,t})  & = \img g_{c,t} + \ker g_{c,t} = \ker h_{c,t} + \img f_{c,t} \\
		\ext^1({\cal T}_r^t, {\cal F}_{m-1}^{c,t}) & = \img f_{c,t} + \ker f_{c,t} = \img f_{c,t} + \img d_{c,t} \\
		\ext^1 ({\cal T}_r^t, {\cal T}_m^t)        & = \ker h_{c,t} + \img h_{c,t}
	\end{align*}

	from which $\ext^1 ({\cal T}_r^t, {\cal F}_{m}^{c,t}) = \ext^1 ({\cal T}_r^t, {\cal T}_m^t) - \img h_{c,t} + \ext^1({\cal T}_r^t, {\cal F}_{m-1}^{c,t}) - \img d_{c,t}$. Both $h$ and $d$ are maps between vector bundles, since by inductive hypothesis the fibers do not change in dimension, and both maps depend only on $\delta$: $d=\delta\circ$ and $h=\delta[1]\circ$. By construction of $\delta$ as a constant extension, the matrix that represents these maps does not depend on $t,c$. Therefore the dimensions of $\img d_{c,t}, \img h_{c,t}$ do not depend on $t,c$ and therefore neither does the dimension of $\ext^1 ({\cal T}_r^t, {\cal F}_{m}^{c,t})$.
\end{proof}
\subsection{Extension to deformation type}\label{sec:ext}
The extension to the entire deformation type follows the same argument as in \cite[Section 4]{Bakker}.
For this reason, we will just sketch the main ideas.\\

Take an \hk manifold $M$ of \Kumm type. Recall from \cref{MukaiHom} that we have an embedding
\begin{equation*}
	H^2(M,\Z)\hookrightarrow\tilde{\Lambda}(M)
\end{equation*}
along with its dual $\dual\theta:\tilde{\Lambda}\rightarrow H_2(M,\Z)$.
Let $\vv(M)$ denote a primitive generator of $H^2(M,\Z)^\perp\subset\tilde{\Lambda}(M)$ (which is of rank 1, as noted in \cite[Rmk. 4.2]{Wieneck}).

An explicit description of the Mori cone is provided by:
\begin{thm}[{\cite[Theorem 2.9]{Mong}}] \label{descrMori}
	Let $(M,h)$ be a polarised holomorphic symplectic variety of generalised Kummer type. The Mori cone of $M$ is generated by the positive cone and classes of the form
	\begin{equation*}
		\left\{\dual{\theta}(a)| a\in\tilde{\Lambda}(M)_{alg }, a^2\geq 0, |(a,\vv(M))|\leq \frac{\vv(M)^2}{2}, h.\theta(a)>0\right\}
	\end{equation*}
\end{thm}
It follows that
\begin{thm}
	Let $(M,h)$ be a holomorphic symplectic variety of \Kumm type and dimension $2n$, with $\PP$ a Lagrangian plane, and let $R\in H_2(M,\Z)$ be its class of line. Then $\tilde\Lambda(M)$ admits a $\PP$-type sublattice $\cal H$ and $R=\dual\theta(s)$ with $a\in\cal H$ a \pit class with $|(v,a)|=\frac{v^2}{2}$
	\proof
	Following the argument in \cite[Theorem 22]{Bakker}, we can deform $M$ to a moduli space, deforming $\PP$ to an extremal Lagrangian plane. On such moduli spaces, we have the corresponding $\PP$-type sublattice $\Hcal$ of the Mukai lattice by \cref{charLag}. By parallel transport, we obtain the $\PP$-type sublattice in $\tilde{\Lambda}(M)$. The same applies for the class of the line $R$.\endproof
\end{thm}
Recall that the discriminant group of a lattice $\Lambda$ is the quotient $\faktor{\Lambda^*}{\Lambda}$ where $\Lambda^*:=\Hom(\Lambda,\Z)$ is the dual lattice.
We can compute the following:
\begin{lemma}
	If $\tilde{\Lambda}(M)$ admits a \pit class $a\in \tilde{\Lambda}(M)$ such that $(a,v)=\frac{v^2}{2}$, then $R=\dual{\theta}(a)$ has $(R,R)=-\frac{n+1}{2}$ and order 2 in the discriminant group of $H^2(M,\Z)$.
	\proof
	By the definition of Mukai homomorphism, we have
	\begin{equation*}
		(R,R)=(a-\frac{v}{2})^2=a^2-(a,v)+\frac{v^2}{4}=-\frac{v^2}{4}
	\end{equation*}
	We also know that:
	$			\frac{v^2}{2}=n+1$
	so we get:
	\begin{equation*}
		(R,R)=-\frac{v^2}{4}=-\frac{n+1}{2}
	\end{equation*}
	Clearly, $2a-v\in v^\perp$ so $R$ is 2-torsion ($2R\in H^2(X,\Z)=v^\perp \subset\tilde{\Lambda}(M)$), but $R\neq 0$ in the discriminant group, so $R$ has order 2.
	\endproof
\end{lemma}
\begin{cor}\label{lineNum}
	Let $M$ be a holomorphic symplectic manifold of \Kumm type and dimension $2n$, and let $R$ be the class of a line of a Lagrangian plane. Then $(R,R)=-\frac{n+1}{2}$ and $2R\in H^2(M,\Z)$.
\end{cor}
\begin{rmk}
	Our result is in alignment with that of Hasset and Tschinkel  \cite{Kumm4fold} on Kummer fourfolds: in that case $n=2$ and so $(R,R)=-\frac{3}{2}$.
\end{rmk}
Another consequence is the following, corresponding to \cite[Corollary 26]{Bakker}
\begin{cor}\label{uniqueMonOrbit}
	There is a single monodromy orbit containing all primitive classes
	arising as lines in Lagrangian planes embedded in holomorphic symplectic varieties of generalised Kummer type.
	\proof The same argument of Bakker holds, as the isometries involved in Eichler's criterion \cite[Chapter 10]{Eichler} have determinant 1 and are orientation preserving, so from \cite[Theorem 1.4]{MarkmannMonKummer} they are monodromy operators also in the generalised Kummer deformation type.\endproof
\end{cor}
The partial converse is proved similarly to \cite[Theorem 25]{Bakker}:
\begin{thm} \label{final}
	Let $(M,h)$ be a polarised holomorphic symplectic manifold of Kummer type and dimension $2n$, with $R \in H_2(M,\Z)$ a primitive generator of an extremal ray of the Mori cone. Then $R$ is the class of a line in a Lagrangian plane if and only if $(R,R) = -\frac{n+1}{2}$ and $2R \in H^2(M,\Z)$.
\end{thm}
\begin{proof}
	One direction of the implication is a direct consequence of the previous theorem.

	The extremality of $R$ implies that, due to the description of the Mori cone, $R=m\dual \theta(a)$ for some $a$ as in \cref{descrMori}, and primitivity forces $m=1$.

	To construct the $\PP$-type lattice $\Hcal$, we use \cref{bastaLaScomposizione}: defining $b:=v-a$ gives the decomposition that proves $\Hcal:=\operatorname{Span}_\Z(a,b)$ is of $\PP$-type. The only thing to verify is that $(a,v)=\frac{v^2}{2}$, which follows from the numerical condition:
	\begin{equation*}
		-\frac{n+1}{2}=(R,R)=(a-\frac{(a,v)}{v^2}v)^2=-\frac{(a,v)^2}{v^2}
	\end{equation*}
	which implies:
	\begin{equation*}
		\frac{(a,v)^2}{v^2}=\frac{n+1}{2}=\frac{v^2}{4}
	\end{equation*}
	so $(a,v)=\frac{v^2}{2}$, and therefore $\Hcal$ is of $\PP$-type.

	Now, as in the previous theorem, we have a smooth proper family along which $R$ stays algebraic that specialises to a moduli space where the image of $R$ is extremal. Thus, we have an extremal Lagrangian plane on the moduli space. Deforming it back to $M$ gives the desired Lagrangian plane $\PP$ having class of a line $R$: as in \cite{Bakker}, since $R$ is primitive, the deformation of the Lagrangian plane to $\PP$ cannot degenerate to something reducible.
\end{proof}
\subsection{Number of planes}
In \cref{sec:construct} we have produced a disjoint union of Lagrangian planes sharing the same class of a line. It's natural to ask about their number. A similar computation has already been performed in an example from \cite[Sect. 4.4]{YoshMod}. We perform Yoshioka's computations with methods and language closer to the ones of this work. We show that, if the Lagrangian planes are extremal, we can reduce to the same setting.

In the following, on an abelian surface $S$ we write $v=(r,D,s), r,s\in\Z,D\in H_\Z^{1,1}(S)$ for the degree decomposition of a Mukai vector $v\in \Hmuk(S,\Z)\cong \Z\oplus H_\Z^{1,1}(S)\oplus \Z$.
\begin{lemma}[Yoshioka's example]\label{exampleYoshioka}
	Let $S$ be an abelian surface with a primitive divisor $D\in \op{Pic}(S), D^2=0$. Consider the Mukai vector $v:=(n+1,D,-1)\in \Hmuk(S,\Z)$ and let $\sigma\in\Stab(S)$ be a $v$-generic stability condition. The Lagrangian planes associated to the partition $[v=(n+1,D,0)+(0,0,-1)]$ as from \cref{construction} are $(n+1)^2$ Lagrangian planes sharing the same class of a line as $\PP$.
	\proof
	Denote $a:=(0,0,-1)$ and $b:=(n+1,D,0)$. Notice that $a^2=b^2=0$ and $(a,b)=n+1$, so it is in the same form as partitions appearing in \cref{construction}. A section of the $\PP^n$ bundle associated with this decomposition in $\M{\sigma}{v}$ is given by the direct sum $\gamma:\M{\sigma}{a}\times\M{\sigma}{b}\rightarrow\M{\sigma}{v}$. Therefore, to count the number of fibres contained in $\K{\sigma}{v}$ one can compute the degree of $\alpha\circ\gamma:\M{\sigma}{a}\times\M{\sigma}{b}\rightarrow S\times \dual{S}$, where $\alpha$ is the Albanese map of $\M{\sigma}{v}$. Notice that this is the Albanese map of the product (as the Albanese map of a product is the sum of the Albanese maps). Moreover, one can observe that $\M{\sigma}{a}\cong \dual S$. Under this identification, we can write the Albanese map explicitly as
	\begin{align*}
		\dual S\times \M \sigma {((n+1,D,0))} & \rightarrow S\times\dual S                                                 \\
		(x, E)                                & \mapsto  \det E-\det E_0,x+\det \Phi_{\cal P}(E)-\det \Phi_{\cal P} (E_0).
	\end{align*}
	To compute the degree it is enough to count the number of solutions to the equations \begin{align*}
		\det E & =\det E_0                                   \\
		x      & =-\det \Phi_{\cal P}(E)+\det \Phi_{\cal P}.
	\end{align*}
	which is the same as the number of solutions of the first equation. This has been computed in \cite{MukaiSemiHom} and it is $(n+1)^2$. \endproof
\end{lemma}
\begin{thm}
	Let $K$ be a projective manifold of generalised Kummer type of dimension $2n$ with an extremal Lagrangian plane $\PP$. Assume that the class $R\in H_2(K,\Z)$ of a line contained in $\PP$ is primitive. Then the number of Lagrangian planes with class of a line $R$ is $(n+1)^2$.
	\proof
	Let $l=\dual{R}\in H^2(K,\Z)$ be the dual of the class $R$.
	We will denote $K_Y$ the moduli space as in \cref{exampleYoshioka}, $\PP_Y$ the corresponding Lagrangian plane and $R_Y\in H^2(K_Y,\Z)$ the class of a line in $\PP_Y$. Denote $l_Y=\dual R_Y$
	By \ref{uniqueMonOrbit}, there is a parallel transport operator $T:H^2(K,\Z)\rightarrow H^2(K_Y,\Z)$ such that $T(l)=l_Y$.
	The operator $T$ can be realised by a path $\gamma:[0,1]\rightarrow {\frak{M}}_n$ in the moduli space of marked IHS manifolds of generalised Kummer type of dimension $2n$. Let $\Lambda$ be the abstract BBF lattice of a manifold of generalised Kummer type. For $t\in [0,1]$, we denote $T_t:\Lambda\rightarrow\Lambda$ the parallel transport operator associated to $\gamma_{|[0,t]}$. In what follows we omit the marking from the notation. Since $T_0(l)=l,T_1(l)= l_Y$ are algebraic, from the same reasoning as in \cite[Corollary 7.4]{MarkmannSurvey} the path can be assumed to preserve algebricity of the class $T_t(l)$ at every point. Moreover, since the classes $l,\tilde l$ are extremal respectively on $K,\tilde{K}$, we can also assume $g_t(l)$ to be extremal for all $t\in[0,1]$: being extremal is an open condition 
	so the locus in ${\frak{M}}_n$ on which we lose extremality has real codimension at least 2 and can be avoided up to modifications of $\gamma$ preserving its homotopy class.

	Since $T_t(l)$ is an extremal algebraic class of negative square, we have a contraction $K_t\rightarrow\hat K_t$ contracting $T_t(l)$.
	As we have shown in \cref{construction}, the contracted locus of the contraction $K_1\rightarrow \hat K_1$ is the disjoint union of the isolated lagrangian planes sharing class of a line $l_Y$. Therefore, the target of the contraction of these planes on the moduli space has isolated singularities, each one corresponding to one of the $(n+1)^2$ contracted lagrangian planes.
	From \cite[Proposition 2.24]{BakkerLehnTheoryArxiv}, locally, deformations of $K_t$ preserving the algebricity of $T_t(l)$ correspond to locally trivial deformation of the target $\hat K_t$. Therefore locally in $t$ the varieties $\hat K_t$ have the same singularities. As $[0,1]$ is connected, we conclude all $\hat K_t$ are singular in $(n+1)^2$ isolated points, meaning that no lagrangian plane can collide with another. As each lagrangian plane is rigid, the contracted locus is a disjoint union of $(n+1)^2$ isolated lagrangian planes.

	\endproof
\end{thm}

\printbibliography
\end{document}